\documentclass{article}
\usepackage{amsfonts}
\newtheorem{thm}{Theorem}

\begin{document}

\title{A George Szekeres Formula for Restricted Partitions}

\author{L. Bruce Richmond\\
Dept. Comb. and Opt., U. of Waterloo\\
Waterloo, Ontario N2L 3G1}
\date{\today}
\maketitle
\begin{abstract}
We derive an asymptotic formula for $A(n,j,r)$ the number of integer partitions of $n$ into at most $j$ parts each part $\le r$. We assume $j$ and $r$ are near their mean values. We also investigate the second
largest part, the number of parts $\ge 2$, etc. We show that the fraction of the partitions of an even integer $n$ that are graphical, ie. whose parts form the degree sequence of a simple graph, is $O(\ln^{-1/2}
n)$. Probabilistic results are used in our discussion of graphical partitions. The George Szekeres circle method is essential for our asymptotic results on partitions. We determine the distributions defined by the successive ranks of partitions, generalizing the result of Erd\H{o}s and Richmond for the rank of a partition.
\end{abstract}

\section{Background}

We let $A(n,j,r) = A(n,r,j)$ denote the number of integer partitions of $n$ into at most $j$ parts each $\le r$. Note $A(n,n,r) = P(n,r)$ the number of partitions $n$  with each part $\le r$. $P(n,r)$ has been investigated by G. Szekeres in \cite{Sz1} and \cite{Sz2}. We shall use techniques and results of Szekeres found in these papers.  We start from the generating function for $A(n,j,r)$ given by Almkvist
and Andrews \cite{AA} Eq.(3.9), by Andrews \cite{Andrews} p. 33, by Comtet \cite{Comtet} p. 117, Ex. 8, by Goulden and Jackson \cite{Gould}, by Stanley \cite{Stanley} p. 29 and by Tak\'{a}cs\cite{Takacs}
$$\sum_{n = 0}^{\infty}A(n,j,r)x^{n}
= \prod_{\nu = 1}^{r}\frac{1 - x^{j + \nu}}{1 - x^{\nu}}.$$
Note $\prod_{\nu = 1}^{r}(1 - x^{\nu})^{-1}$ is the generating function for $P(n,r)$. The generating function of $A(n,j,r)$ is the classical Gauss binomial coefficient, see Gauss \cite{Gauss}, p. 16.

We use what can be called the Szekeres circle method. This method is an example of the saddlepoint method applied to partitions. One writes $A(n,j,r)$
as a contour integral using Cauchy's theorem
\begin{equation}\label{I}
A(n,j,r) = \frac{1}{2\pi i}\int_{C}\prod_{\nu = 1}^{r}\frac{1 - x^{j + \nu}}{1 - x^{\nu}}x^{-n - 1}dx
\end{equation}
where the contour of integration $C$ is a circle of radius $e^{-\alpha}$, that is $x = e^{-\alpha + i \theta}$ and we shall modify the Szekeres choice of $\alpha$ slightly, we choose 
$$\alpha = \frac{\pi}{\sqrt{6n}}.$$
Szekeres \cite{Sz1} shows that the integral over a small arc $|\theta| \le \theta_{o},\; \theta_{o} = n^{-5/7}$,
see Eq. (15) of \cite{Sz1}, gives the asymptotic behaviour of $P(n,r)$ and this true for our integral also.

The saddlepoint method is discussed in Flajolet and Sedgewick \cite{Flajolet} for example. If the integrand in Eq~(\ref{I}) is expressed as (setting $x = e^{-\alpha + i\theta}$)
$$e^{f(\alpha) - \frac{1}{2}f^{''}(\alpha)\theta^{2} + \frac{1}{6}i^{3}f^{'''}(\alpha)\theta^{3} + \cdots}$$
($\alpha$ is chosen so that the coefficient of $\theta$ is zero) then
\begin{equation}\label{II}
\int_{-\theta_{o}}^{\theta_{o}}e^{f(\alpha + i\theta)}d\theta = \frac{e^{f(\alpha)}}{\sqrt{2\pi f^{''}(\alpha)}}\left(1 + O\left(f^{'''}(\alpha)\theta_{o}^{3}\right)\right).
\end{equation}
For this estimate for $A(n,j,r)$ to hold it is necessary that the integral over those $\theta \in [-\pi, -\theta_{o}] \cup [\theta_{0}, \pi]$ is negligable. The arguments of Szekeres \cite{Sz1} suffice to show that the integral over those $\theta \in [-\pi,-\theta_{0}] \cup [\theta_{0}, \pi]$ is negligable in our case as well. We shall show that the $O$-term in Eq.~(\ref{II}) is $O(n^{-1/7})$.

We shall see that the distributions defined by the largest summand and the number of summands are independent. Szekeres \cite{Sz2} states that this result is not available. Erd\H{o}s and Richmond \cite{Erdos}
verify that this is true
however we find a more careful analysis of the error terms necessary for our present results. 

To apply the saddlepoint method as described in Eq.~(\ref{II}) we consider the function $G(\theta)$ to be the integrand in Eq.~(\ref{I}) then consider the logarithm of $G(\theta)$. We have letting $x = e^{-\alpha + i\theta}$
$$\ln G(e^{-\alpha + i\theta}) = -\sum_{\nu = 1}^{r}\ln\left(1 - e^{(-\alpha + i\theta)\nu}\right) + \sum_{\nu = 1}^{r}\ln\left(1 - e^{(-\alpha + i\theta)(\nu + j)}\right) - i n\theta$$
so the derivative of the logarithm of $G$ at $\theta = 0$ is $i$ times
$$ \sum_{\nu = 1}^{r}\frac{\nu}{e^{\alpha \nu} - 1} - \sum_{\nu = 1}^{r}\frac{\nu + j}{e^{\alpha(\nu + j)} - 1} - n.$$
Thus choosing $\alpha$ so that saying the  coefficient of $\theta$ is 0 is equivalent to choosing $\alpha$ so that 
$$n = \sum_{\nu = 1}^{r}\frac{\nu}{e^{\alpha \nu} - 1} - \sum_{\nu = 1}^{r}\frac{\nu + j}{e^{\alpha(\nu + j)} - 1}.$$
We shall see however that our previous choice of $\alpha$ works for our purposes. The second sum turns out to be negligable compared to the first.

Almkvist and Andrews investigate the asymptotic behaviour of $A(j,n,r)$ when $j$ is near $nr/2$, these values of $j$ are the middle values since $A(j,n,r) = A(nr - j,n,r)$ for $j = 0,1, \cdots, nr$. They point out that $A(j,n,r)$ has applications in statistics, the Wilcoxon rank sum test.

L. Tak\'{a}cs \cite {Takacs} interprets $A(n,j,r)$ as the number of lattice paths from $(0,0)$ to $(n,j)$ the area under which is $r$.
 To verify that the generating function for $A(n,j,r)$ is what we state Tak\'{a}cs derives 
 a recurrence. In terms of partitions his recurrence comes from the observation that with $A(0,0,0) = 1, A(j, 0, 0) = 0$ for $j \ge 2$ the partitions of $n - r$ with largest part $\le r$ and having at most $j - 1$ parts plus the partitions of $n$ having largest part $\le r - 1$ and having at most $j$ parts are in 1-1 correspondence with the partitions counted by $A(n,j,r)$. This is clear from the Ferrar's diagram for a partition. Thus, if $n < 0$ then $A(n,j,r) = 0$,
$$\sum_{n = 0}^{\infty}A(n,j,r)z^{n} = z^{r}\sum_{n = 0}^{\infty}A(n - r,j - 1,r)z^{n - r} + \sum_{n = 1}^{\infty}A(n,j,r-1)z^{n}.$$
If we let
\[
\sum_{0}^{\infty}A(n,j,r)z^{n} = \left[
\begin{array}{c}
j + r\\
j
\end{array}
\right]_{z} =
\left[
\begin{array}{c}
j + r\\ r
\end{array}
\right]_{z}
\]
this recurrence becomes
\[
z^{r}\left[
\begin{array}{c}
j - 1 + r\\j - 1
\end{array}
\right]_{z}
+\;
\left[
\begin{array}{c}
j + r - 1\\j 
\end{array}
\right]_{z},
\]
which is one of the recurrences Gauss gives for his binomial coefficient.

R. Stanley \cite{Stanley} p. 30 interprets $A(n,j,r)$ as the number of permutations of the multiset $\{ 1^{j}, 2^{r}\}$ with $n$ inversions.

 Tak\'{a}cs shows that $A(n,j,r)$ defines a normal distribution when $j$ is near $nr/2$, giving the mean$(=nr/2)$ and standard deviation.
(Almkvist and Andrews also do this by a different method). It is easy to see using the previous product formula for the generating function of $A(n,j,r)$ that that the generating function satisfies this recurrence.

 Since $A(j,n,r) = A(nr - j,n,r)$ the authors just mentioned investigate $A(j,n,r)$ when $j$ is near the midpoint $nr/2$ and show that this midpoint is 
 the mean value of a distribution which tends to a normal distribution. This is of course a natural thing to do. We investigate $A(n,j,r)$ with $j$ and $r$ near $\pi^{-1}\sqrt{6n}\ln(\pi^{-1}\sqrt{6n})$, see Theorem~\ref{one}.
Thus we investigate $A(n,j,r)$ for large $n$ near its maximum value viewed as a function of the two variables $j$ and $r$. $A(n,j,r)$ defines a distribution which is the product of two identical independent extremal distributions, again  see Theorem~\ref{one}.

In Section 3 we generalize the result of Erd\H{o}s and Richmond for the rank of a partition to successive ranks.

\section{Restricted Partitions}
We begin with Eq.~(\ref{II}) and determine $f(\alpha)$.
Consider the $\prod_{\nu = 1}^{r}(1 - e^{(-\alpha + i\theta) \nu})^{-1}$ part of Eq.~(\ref{I}) first. At the end of Section 2 of \cite{Sz1}, p 105, we find with $u = \alpha r$(our $r$ is denoted by $k$ in this paper of Szekeres)
$$-\sum_{\nu = 1}^{r}\ln(1 - e^{-\nu \alpha}) = \frac{1}{\alpha}\int_{0}^{u}\frac{t}{e^{t} - 1}dt -\left( \frac{u}
{\alpha} + \frac{1}{2}\right)\ln\left(1 - e^{-u}\right) + \frac{1}{2}\ln\frac{\alpha}{2\pi} + O(\alpha).$$
To estimate this integral, let's denote it by $I$, we note that using integration by parts we have, since
$$ \int\frac{t}{e^{t} - 1}dt = t\ln\left(1 - e^{-t}\right) + \int \ln\left(\frac{1}{1 - e^{-t}}\right)dt,$$
that
\begin{equation}\label{Y}
I = \int_{0}^{\alpha r}\frac{t}{e^{t} - 1}dt = -\alpha re^{-\alpha r} +O\left(\alpha re^{-2\alpha r}\right) +
 \sum_{l = 1}^{\infty}-\frac{e^{-\alpha r l}}{l^{2}} + \sum_{l = 1}^{\infty}\frac{1}{l^{2}}
\end{equation}
$$ = \frac{\pi^{2}}{6} -\alpha re^{-\alpha r} - e^{-\alpha r} + O\left(\alpha r e^{-2\alpha r}\right)$$
provided $\alpha r \rightarrow \infty$.

 Since it is known that the expected value of $r$ is $\pi^{-1}(6n)^{1/2}\ln(\pi^{-1}\sqrt{6n}) $, see \cite{Sz2}, we let $r = \pi^{-1}\sqrt{6n}\ln(\pi^{-1}\sqrt{6n}) + x_{1}\pi^{-1}\sqrt{6n}$ and we find 
$-\alpha r = -\ln(\pi^{-1}\sqrt{6n}) - x_{1} = O(|x_{1}| + \ln n).$ Thus
$$e^{-\alpha r} = \frac{\pi}{\sqrt{6n}}e^{- x_{1}}.$$
Finally we have the following estimate for our integral $I$, provided that $|x_{1}| < (1 - \epsilon) \ln(\pi^{-1}\sqrt{6n})$,

\begin{equation}\label{III}
I = \frac{\pi^{2}}{6} - \frac{\pi}{\sqrt{6n}}e^{-x_{1}} -
\alpha r \frac{\pi}{\sqrt{6n}}e^{-x_{1}} + O\left(n^{-1}(|x_{1}| + \ln n)e^{- 2x_{1}}\right).
\end{equation}

It follows that
$$\alpha^{-1}\int_{0}^{\alpha r}\frac{t}{e^{t} - 1}dt = \pi \sqrt{n/6} - e^{-x_{1}} - r\frac{\pi}{\sqrt{6n}}e^{-x_{1}} + O(n^{-1/2}(|x_{1}| + \ln n)e^{-2x_{1}}).$$
We now have from the last equation of Section 2 of \cite{Sz1}
$$-\sum_{\nu = 1}^{r}\ln\left(1 - e^{-\nu \alpha}\right) = \frac{1}{\alpha}\int_{0}^{r\alpha}\frac{t}{e^{t} - 1}dt - \left(\frac{r\alpha}{\alpha} + \frac{1}{2}\right)\ln\left(1 - e^{-r \alpha}\right)$$
$$ + \frac{1}{2}\ln \frac{\alpha}{2\pi} + O(\alpha).$$

It follows that, the $r\pi(6n)^{-1/2}e^{-x_{i}}$ term cancells,
\begin{equation}\label{IV}
-\sum_{\nu = 1}^{r}\ln \left(1 - e^{-\nu \alpha}\right) = \pi \sqrt{n/6} - e^{-x_{1}} + 2^{-1}\ln(\alpha/2\pi) + O(n^{-1/2}(|x_{1}| + \ln n)e^{-2x_{1}}).
\end{equation}

 Let us now consider
$$\prod_{\nu = 1}^{r}\left(1 - x^{j + \nu}\right), \;\; x = e^{-\alpha + i\theta}$$
and suppose since the expected value of $j$ equals the expected value of $r$(there is a 1-1 correspondence between partitions and their conjugates)
$$j = \frac{\sqrt{6n}}{\pi}\ln\left(\frac{\sqrt{6n}}{\pi}\right) + y_{1}\frac{\sqrt{6n}}{\pi},\;\; \alpha j = \ln\left(\frac{\sqrt{6n}}{\pi}\right) + y_{1}.$$
We have
$$\prod_{\nu = 1}^{r}\left(1 - e^{-(\nu + j) \alpha}\right) = e^{\sum_{\nu = 1}^{r}\ln\left(1 - e^{-\alpha j + \alpha \nu}\right)}$$
so we consider, much as before $e^{-\alpha j} = e^{-y_{1}}\pi(6n)^{-1/2}$,
$$\sum_{\nu = 1}^{r}\ln\left(1 - e^{-\alpha j + \alpha \nu}\right) = -\sum_{\nu = 1}^{r} e^{-\alpha j - \alpha \nu} - \frac{e^{-2\alpha(\nu + j)}}{2} - \cdots$$
$$= -e^{-\alpha j}\sum_{\nu = 1}^{r}e^{-\alpha \nu} + O\left(e^{-2\alpha j}\sum_{\nu = 1}^{r}e^{-2\alpha \nu}\right)$$
$$ = -e^{-\alpha j}\frac{1 - e^{-\alpha(r + 1)}}{1 - e^{-\alpha}} + O\left(e^{-2\alpha j}\frac{1 - e^{-2\alpha(r + 1)}}{1 - e^{-2\alpha}}\right).$$ 
Now 
$$\frac{e^{-\alpha j}}{1 - e^{-\alpha}} = e^{-y_{1}}\left(1 + O\left(n^{-1/2}\right)\right)$$
and
$$1 - e^{-\alpha(r + 1)} = 1 + O\left(\left(e^{-x_{1}} + 1\right)n^{-1/2}\right)$$
so we conclude that
\begin{equation}\label{V}
\sum_{\nu = 1}^{r}\ln\left(1 - e^{-\alpha j + \alpha \nu}\right) = -e^{-y_{1}} + O\left(\frac{e^{-x_{1} - y_{1}} + e^{-y_{1}} + 1}{n^{1/2}}\right).
\end{equation}
Thus from Eqs.~(\ref{IV}) and ~(\ref{V})
we have that the $f(\alpha)$ in Eq.~(\ref{II}) satisfies(notice the term $e^{-\alpha n}$ in Eq.~(\ref{II}) )
\begin{equation}\label{VI}
f(\alpha) = \pi\sqrt{2n/3} - e^{-x_{1}} - e^{-y_{1}} + 2^{-1}\ln\left(2^{-3/2}3^{-1/2}n^{-1/2}\right) + O\left(\frac{e^{-x_{1} - y_{1}} + e^{-y_{1}} + 1}{n^{1/2}}\right)
\end{equation}
$$ + O\left(\frac{(|x_{1}| + \ln n)e^{-2x_{1}}}{n^{1/2}}\right).$$

Let us return to the $G(\theta)$ we considered in the background Section. We have that the second derivative of the logarithm of $G(\theta)$ equals
$$ -\sum_{\nu = 1}^{r}\frac{\nu^{2}e^{(\alpha + i\theta)\nu}}{\left(e^{(\alpha + i\theta)\nu} - 1\right)^{2}}
 + \sum_{\nu = 1}^{r}\frac{(j + \nu)^{2}e^{(\alpha + i\theta)
(j + \nu)}}{\left(e^{(\alpha + i\theta)(j + \nu)} - 1\right)^{2}}$$
evaluated at $\theta = 0$ this gives
$$-\sum_{\nu = 1}^{r}\frac{\nu^{2}e^{\alpha \nu}}{\left(e^{\alpha \nu} - 1\right)^{2}} + \sum_{\nu = 1}^{r}\frac{(j + \nu)^{2}e^{\alpha(j + \nu)}}{\left(e^{\alpha (j + \nu)} - 1\right)^{2}}.$$
As before it follows that $\alpha(j + r) = 2\ln(\sqrt{6n}/\pi) + x_{1} + y_{1}$
$$e^{\alpha(j + r)} = \frac{6n}{\pi^{2}}e^{x_{i} + y_{i}}.$$
Note
$$(j + r)^{2}r = \left(\frac{\sqrt{6n}}{\pi}\right)^{3}
\left(4\ln^{2}\frac{\sqrt{6n}}{\pi} + (4x_{1} + 4y_{1})\ln \frac{\sqrt{6n}}{\pi} + (x_{1} + y_{1})^{2}\right)\left(\ln\left(\frac{\sqrt{6n}}{\pi}\right) + x_{1}\right).$$
 Thus if $|x_{1}|, |y_{1}| \le 5^{-1}\ln n$(this implies that our previous bound $|x_{1}| < (1 - \epsilon)\ln(\pi^{-1}\sqrt{6n})$ holds) then, since $e^{\alpha(j + \nu)}/(e^{\alpha(j + \nu)} - 1)^{2} \sim e^{-\alpha(j + \nu)} \le e^{-\alpha j}$ = $O(e^{-y_{1}}n^{-1/2}) = O(n^{-3/10})$
and $|e^{-x_{1} - y_{1}}| \le n^{2/5}$, we have
\begin{equation}\label{VII}
\sum_{\nu = 1}^{r}\frac{(j + \nu)^{2}e^{\alpha(j + \nu)}}{\left(e^{\alpha(j + \nu)} - 1\right)^{2}} = O\left((j + r)^{2}r e^{-(\alpha + j)}\right)
 = O\left(n^{6/5}\ln^{3}n\right).
\end{equation}
From Eq. (12) 
 of \cite{Sz1} we have first
$$\sum_{\nu = 1}^{r}\frac{\nu^{2}e^{\alpha \nu}}{\left(e^{\alpha \nu} - 1\right)^{2}} = \alpha^{-3}\int_{0}^{\alpha r}\frac{t^{2}e^{t}}{\left(e^{t} - 1\right)^{2}}dt + O\left(\alpha^{-2}\right)$$
then Eq.(3) of \cite{Sz1} implies that this equals, $\alpha^{2}r^{2}(e^{\alpha r} - 1)^{-1} = O(\ln^{2}n e^{-\alpha r}) = O(n^{-3/10}\ln^{2}n)$,
$$ = \alpha^{-3}\left[2\int_{0}^{\alpha r}\frac{t}{e^{t} - 1}dt - \frac{\alpha^{2}r^{2}}{e^{\alpha r} - 1}\right] + O\left(\alpha^{-2}\right).$$

Eq.~(\ref{VII}) shows that the second sum in our formula for $(\ln G(\theta))^{''}$
 when $\theta = 0$ is negligable compared to the first which exceeds a constant times $n^{3/2}$. It is also true that the second sum in our formula for $(\ln G(\theta))^{'}(\theta)$ is negligable compared to the first as we stated in the Background section by a similar argument.

Furthermore our Eq.~(\ref{III}) now gives that this last expression is
$$= \left(\frac{\pi}{\sqrt{6n}}\right)^{-3}\left[\frac{\pi^{2}}{3} + O\left(\frac{(\ln n + |x_{1}|)e^{-x_{1}}}{n^{1/2}}\right)\right] = \frac{2\sqrt{6}}{\pi}n^{3/2} + O\left(n(\ln n + |x_{1}|)e^{-x_{1}}\right).$$
This result with Eq.~(\ref{VII}) allows us to estimate $f^{''}(\alpha)$ which appears in Eq.~(\ref{II}). It implies that, the first $O$-term comes from Eq.~(\ref{VII}); 
the second $O$-term is only as big as $O(n^{-3/10}\ln n)$ for $|x_{1}| \le 5^{-1}\ln n$,
$$(2\pi f^{''}(\alpha))^{-1/2} = \frac{1}{6^{1/4}2 n^{3/4}}\left(1 + O\left(\frac{\ln^{3}n}{n^{1/5}}
\right) + O\left(\frac{(\ln n + |x_{1}|)e^{-x_{1}}}{n^{1/2}}\right)\right).$$
Furthermore from Eq.~(\ref{VI}), the $O(n^{-1/5}\ln^{3} n)$ term is negligable,
$$e^{f(\alpha)} = e^{\pi\sqrt{2n/3}}e^{-e^{-x_{1}}}e^{-e^{-y_{1}}}\frac{1}{2^{3/4}3^{1/4}n^{1/4}}\left(1 + O\left(\frac{e^{-x_{1} - y_{1}} + e^{-y_{1}} + (\ln n + |x_{1}|)e^{-2x_{1}}}{n^{1/2}}\right)\right).$$
Hence  for $|x_{1}|, |y_{1}| \le 5^{-1}\ln n$(again the $O(n^{-1/5}\ln^{3}n)$ term is negligable)
\begin{equation}\label{Z}
\frac{e^{f(\alpha)}}{(2\pi f^{''}(\alpha))^{1/2}} = \frac{e^{\pi\sqrt{2n/3}}}{4\sqrt{3}n}e^{-e^{-x_{1}}}
e^{-e^{-y_{1}}}\left(1 + O\left(\frac{e^{-x_{1} - y_{1}} +e^{-y_{1}} + (\ln n + |x_{1}|)e^{-2x_{1}}}{n^{1/2}}\right)\right)
\end{equation}
$$= \frac{e^{\pi\sqrt{2n/3}}}{n4\sqrt{3}}e^{-e^{-x_{1}}}e^{-e^{-y_{1}}}\left(1 + O\left(n^{-1/10}\ln n\right)\right).$$

We have calculated the second derivative of the logarithm of $G(\theta)$ and estimated it at $\theta = 0$ to apply Eq.~(\ref{II}). We also need to bound the third derivative. We find that
$$\left(\ln G(\theta)\right)^{'''} = -i\sum_{\nu = 1}^{r}\frac{\nu^{3}e^{(\alpha + i\theta)\nu}}{\left(e^{(\alpha + i\theta)\nu} - 1\right)^{2}} + 2i\frac{\nu^{3}e^{2(\alpha + i\theta)\nu}}{\left(e^{(\alpha + i\theta)\nu} - 1\right)^{3}}$$
$$ +i\sum_{\nu = 1}^{r}\frac{(\nu + j)^{3}e^{(\alpha + i\theta)(\nu + j)}}
{\left(e^{(\alpha + i\theta)(\nu + j)} - 1\right)^{2}} -2i\sum_{\nu = 1}\frac{(j + \nu)^{3}e^{2(\alpha + i\theta)(j + \nu)}}{\left(e^{(\alpha + i\theta)(\nu + j)} - 1\right)^{3}}.$$
We wish to estimate this when $\theta = 0$. From Eq. (12) of \cite{Sz1} we have
$$\sum_{\nu = 1}^{k}\nu^{r}e^{s\alpha \nu}(e^{\alpha \nu} - 1)^{-r} = \alpha^{- r - 1}\int_{0}^{\alpha k} t^{r}e^{st}(e^{t} - 1)^{-r} dt + O(\alpha^{-r}).$$
We may apply this formula of Szekeres, his $k$ is our $r$ and his $r$ is 3 in our case,  after bringing the first two terms in our formula for the third derivative of $\ln G(\theta)$ to a common denominator and 
recalling that since $\alpha r = \ln(\sqrt{6n}/\pi) \rightarrow \infty$ we find that the first two terms in this formula are $O\left(\alpha^{-4}\right)$. The two terms involving $\nu + j$ are $O((j + r)^{3}re^{-\alpha r}) = O(n^{17/10}\ln^{4}n)$
 for $|x_{1}|,|y_{1}| \le 5^{-1}\ln n$ by the argument giving Eq.~(\ref{VII}) and we conclude that the third derivative of $\ln G(\theta)$ is $O(n^{2})$. Since $\theta_{0} = n^{-5/7}$ we conclude that the $O$-term in our Eq
~(\ref{II}) is $O(n^{-1/7})$.

In view of the famous Hardy-Ramanujan formula for $P(n)$(a weak consequence thereof rather), the number of partitions of $n$, we have proved using Eqs.~(\ref{Z}) and ~(\ref{II}) apart from proving that the integral of $G(\theta)$ over $\theta \in [-\pi, -\theta_{0}]\cup[\theta_{0}, \pi]$ is negligable(which follows easily from results of Szekeres as we shall see), we use that $e^{-2x_{1}}, e^{-x_{1} - y_{1}} = O(n^{2/5})$,

\begin{thm}\label{one}
If $A(n,j,r)$ denotes the number of integer partitions of $n$ into at most $j$ parts each part $\le r$ and 
$$ j = \frac{\sqrt{6n}}{\pi}\ln\left(\frac{\sqrt{6n}}{\pi}\right) + y_{1}\frac{\sqrt{6n}}{\pi}, \;\; r = \frac{\sqrt{6n}}{\pi}\ln\left(\frac{\sqrt{6n}}{\pi}\right) + x_{1}\frac{\sqrt{6n}}{\pi}$$
then if $|x_{1}|, |y_{1}| \le 5^{-1}\ln n$
$$A(n,j,r) = P(n)e^{-e^{-x_{1}}}e^{-e^{-y_{1}}}\left(1 + O\left(\frac{e^{-x_{1} - y_{1}} + e^{-y_{1}} + (\ln n + |x_{1}|)e^{-2x_{1}}}{n^{1/2}}\right)\right)$$
$$= P(n)e^{-e^{-x_{1}}}e^{-e^{-y_{1}}}
\left(1 + O\left(n^{-1/10}\ln n\right)\right).$$
The distributions defined by $j$ and $r$ are thus independent and equal.
\end{thm}

This agrees with the case $d = 1$ of Fristedt's Theorem 6.3 \cite{Fristedt} for the distribution of $r$. The independence result seems to be new.

Let us now show that the integral of $G(\theta)$ over $\theta \in [-\pi, -\theta_{0}]\cup[\theta_{0}, \pi]$ is negligable. Consider the logarithm of $\prod_{\nu = 1}^{r}(1 - x^{\nu})$
 which is part of the Almkvist-Andrews generating function for $A(n,j,r)$, that is
$$\sum_{\nu = 1}^{r}\ln\left(1 - e^{(- \alpha + i\theta)(\nu + j)}\right) = -\sum_{\nu = 1}^{r}\ln\left(\frac{1}{1 - e^{(-\alpha + i\theta)(\nu + j)}}\right)$$
$$ = -\sum_{\nu = 1}^{r}e^{(-\alpha + i\theta)(\nu + j)} - \sum_{\nu = 1}^{r}\frac{e^{2(- \alpha + i\theta)(\nu + j)}}{2} - \cdots.$$
We have 
$$\sum_{\nu = 1}^{r}e^{(-\alpha + i\theta)(\nu + j)} = e^{(-\alpha + i\theta)j}
\frac{1 - e^{(- \alpha + i\theta)(r + 1)}}{1 - e^{-\alpha + i\theta}}.$$
Thus
$$\left|\sum_{\nu = 1}^{r}e^{(-\alpha + i\theta)(\nu + j)}\right| \le e^{-\alpha j}\frac{1 + e^{-\alpha(r + 1)}}{\left|1 - e^{-\alpha + i\theta}\right|}.$$
Moreover

$$e^{-\alpha j} = \frac{\pi e^{-y_{1}}}{\sqrt{6n}}$$
so if $|y_{1}| \le 5^{-1}\ln n$ then $e^{-\alpha j} = O(n^{-3/10}).$  We also have
$$\left|\frac{1}{1 - e^{-\alpha + i\theta}}\right| \le \frac{1}{1 - e^{-\alpha}\cos \theta}
\le \frac{1}{1 - e^{-\alpha}\cos \theta_{0}}.$$
Now
$$e^{-\alpha} = 1 -\frac{\pi}{\sqrt{6n}} + O(n^{-1}), \;\; \cos \theta_{0} = 1 - n^{-10/7}/2 + O(n^{-20/7})$$
so
$$e^{-\alpha}\cos \theta_{0} = 1 -\frac{\pi}{\sqrt{6n}} + O(n^{-1}),
 \;\; 1 - e^{-\alpha}\cos \theta_{0} = \frac{\pi}{\sqrt{6n}} + O(n^{-1})$$
and
$$\left|\frac{1}{1 - e^{-\alpha + i\theta}}\right| \le \frac{\sqrt{6n}}{\pi} + O(1) = \frac{\sqrt{6n}}{\pi}\left(1 + O\left(n^{-1/2}\right)\right).$$
We now combine these facts, since $e^{-\alpha(r + 1)} = \pi e^{-x_{1}}/\sqrt{6n}(1 + O(n^{-1/2}))$,
$$e^{-\alpha j}\frac{1 + e^{-\alpha(r + 1)}}{\left|1 - e^{-\alpha + i\theta}\right|} \le \frac{\pi}{\sqrt{6n}}e^{-y_{1}}
\frac{\sqrt{6n}}{\pi}(1 + O(n^{-1/2}))(1 + \pi/\sqrt{6n}e^{-x_{1}}(1 + O(n^{-1/2})))$$
$$ = \left(e^{-y_{1}} + \frac{\pi}{\sqrt{6n}}e^{-x_{1} -y_{1}}\right)\left(1 + O\left(n^{-1/2}\right)\right)$$
which is $O(n^{1/5})$ if $|y_{1}| \le 5^{-1}\ln n$.
The other sums in the logarithm of this product can be estimated in a very similar way. The $l$-th term in our sum for $\ln(1 - e^{(-\alpha + i\theta)(\nu + j)})$ is
$$l^{-1}e^{(-\alpha + i\theta)lj}\frac{1 - e^{(- \alpha + i\theta)l(r + 1)}}{1 - e^{l(-\alpha + i\theta)}}.$$
Now $e^{-l\alpha j} = e^{-l \ln(\sqrt{6n}/\pi) - ly_{1}} = (\pi e^{-y_{1}}/\sqrt{6n})^{l}$. 
 Furthermore
$$\left|\frac{1}{1 - e^{(-\alpha  + i\theta)l}}\right|^{2} = \frac{1}{\left(1 -e^{-\alpha l}\right)^{2} + 2e^{-\alpha l}(1 - \cos l\theta)} \le \left(\frac{1}{1 - e^{-\alpha l}}\right)^{2}.$$
This gives a bound for the $l$-th sum of
$$\left(\frac{\pi e^{-y_{1}}}{\sqrt{6n}}\right)^{l}\frac{1}{1 - e^{-\alpha}} = \left(\frac{\pi e^{-y_{1}}}
{\sqrt{6n}}\right)^{l}\left(\frac{\sqrt{6n}}{\pi} + O(1)\right).$$
If $|y_{1}|\le 5^{-1}\ln n$  this bound is $O(n^{-3l/10 + 1/2})$ and so the sum over all terms $l \ge 2$ is $O(n^{-1/10})$. Thus the sum over all terms is $O(n^{1/5})$.
 In other words $|\prod_{\nu = 1}^{r}(1 - x^{\nu})| = O(n^{1/5})$ on the path of integration $|x| = e^{-\alpha}.$ Now Szekeres \cite{Sz1} shows, Eq.(23), that for $|\theta| \ge \theta_{0} = n^{-5/7}$ the absolute value of $\prod_{\nu = 1}^{r}(1 - x^{\nu})^{-1}$ is bounded by $e^{-cn^{1/14}}$ times the value of this product evaluated at $|x| = e^{-\alpha}$.
 The same bound holds for the Almkvist-Andrews generating function for $A(n,j,r)$ with
a smaller constant $c_{1}$ than the Szekeres constant, that is if $|\theta| \ge \theta_{0}$ then $|G(\theta)| \le G(0)e^{-c_{1}n^{1/14}}$. Thus all the conditions for Eq.~(\ref{II}) are satisfied and it holds, 
the proof of Theorem~\ref{one} is complete.

 We now show that we may modify
the restriction that $|x_{1}|, |y_{1}| \le 5^{-1}\ln n$ in this theorem. Note if $x_{1} = 5^{-1}\ln n$
then $e^{-e^{-x_{1}}} = e^{-n^{-1/5}} = 1 + O(n^{-1/5})$. If $N(x_{1})$ denotes the number of partitions of $n$ with largest part $\le
\frac{\sqrt{6n}}{\pi}\ln(\frac{\sqrt{6n}}{\pi}) + x_{1}\frac{\sqrt{6n}}{\pi}$ with $
x_{1} \ge 5^{-1}\ln n$ then since $N(x_{1})$ increases with $x_{1}$ we have from Theorem ~\ref{one}
$$ P(n)(1 + O(n^{-1/5})) \le N(x_{1}) \le P(n)$$ 
so Theorem ~\ref{one} holds for $x_{1} \ge 5^{-1/5}\ln n$(since we may replace $e^{-e^{-x}}$ by 1 for $x_{1} \ge 5^{-1}\ln n)$. A similar argument allows $y_{1} \ge 5^{-1}\ln n$. 
Finally if $x_{1} = -5^{-1}\ln n$ then $e^{-e^{-x_{1}}} = e^{-n^{1/5}}$ and 
if $N1(x_{1})$ denotes the number of partitions of $n$ with largest part $\le \frac{\sqrt{6n}}{\pi}\ln(\frac{\sqrt{6n}}{\pi}) + 
x_{1}\frac{\sqrt{6n}}{\pi}$ we find that if $x_{1} \le -5^{-1}\ln n$ then $N1(x_{1}) \le e^{-n^{1/5}}$. The same bound holds for $y_{1} \le -5^{-1}\ln n$.

We now have
\begin{thm}\label{two}
If $A(n,j,r)$ denotes the number of integer partitions of $n$ into at most $j$ parts each part $\le r$ and
$$j = \frac{\sqrt{6n}}{\pi}\ln\left(\frac{\sqrt{6n}}{\pi}\right) + y_{1}\frac{\sqrt{6n}}{\pi}, \;\; r = \frac{\sqrt{6n}}{\pi}\ln\left(\frac{\sqrt{6n}}{\pi}\right) + x_{1}\frac{\sqrt{6n}}{\pi}$$
then if $x_{1} \ge -5^{-1}\ln n$ and $y_{1} \ge -5^{-1}\ln n$ then 
$$A(n,j,r) = P(n)e^{-e^{-x_{1}}}e^{-e^{-y_{1}}}\left(1 + O\left(n^{-1/10}\ln n\right)\right).$$
If one of $x_{1}$ or $y_{1}$ is $\le -5^{-1}\ln n$ then
$$A(n,j,r) = O\left(P(n)e^{-n^{1/5}}\right).$$
\end{thm}
 
 We now study the number of partitions with largest part $r$ and having $k$ parts. Let $B(n,k,r) = A(n,k,r) - A(n,k,r - 1)$ be the number of partitions of $n$ with largest part $r$ and number of parts $\le k$.
From the Almkvist-Andrews generating function we have
$$\sum_{n = 0}^{\infty}B(n,k,r)x^{n} = \frac{1 - x^{k + r}}{1 - x^{r}}
\prod_{\nu = 1}^{r - 1}\frac{1 - x^{k + \nu}}{1 - x^{\nu}} - \prod_{\nu = 1}^{r - 1}\frac{1 - x^{k + \nu}}{1 - x^{\nu}}$$
$$= \left(\frac{1 - x^{k + r}}{1 - x^{r}} - 1\right)\prod_{\nu = 1}^{r - 1}
\frac{1 - x^{k + \nu}}{1 - x^{\nu}} = \frac{x^{r}}{1 - x^{r}}\left(1 - x^{k}\right)\prod_{\nu = 1}^{r - 1}\frac{1 - x^{k+ \nu}}{1 - x^{\nu}}.$$

If we temporarily let $O(x^{k})$ denote a formal power series that has smallest power of $x$ equal to $k$ we may now write
$$\sum_{n = 0}^{\infty}B(n,k,r)x^{n} = \left(x^{r} + x^{2r} - x^{k + r} + O(x^{3r}) + O(x^{k + 2r})\right)\prod_{\nu = 1}^{r - 1}\frac{1 - x^{k + \nu}}{1 - x^{\nu}}.$$

Now let $C(n, k, r) = B(n,k,r) - B(n, k - 1,r)$ denote the number of partitions of $n$ with largest part = $r$ and having exactly $k$ parts. Then
$$C(x) = \sum_{n = 0}^{\infty}C(n,k,r)x^{n} = \frac{x^{r}}{1 - x^{r}}\left(1 - x^{k}\right)\prod_{\nu = 1}^{r - 1}\frac{1 - x^{k + \nu}}{1 - x^{\nu}}$$
$$ - \frac{x^{r}}{1 - x^{r}}\left(1 - x^{k - 1}\right)\prod_{\nu = 1}^{r - 1}\frac{1 - x^{k - 1 + \nu}}{1 - x^{\nu}}.$$
So
$$C(x) =\left[ \frac{x^{r}}{1 - x^{r}}\left(1 - x^{k}
\right) - \frac{x^{r}}{1 - x^{r}}\frac{(1 - x^{k - 1})(1 - x^{k})}{1 - x^{k + r - 1}}\right]\prod_{\nu = 1}^{r - 1}\frac{1 - x^{k + \nu}}{1 - x^{\nu}}$$
$$=\frac{x^{r}}{1 - x^{r}}\left(1 - x^{k}\right)\left[\frac{x^{k - 1} - x^{k + r - 1}}{1 - x^{k + r - 1}}\right]\prod_{\nu = 1}^{r - 1}\frac{1 - x^{k + \nu}}{1 - x^{\nu}}$$
$$= x^{r + k - 1}\frac{1 - x^{k}}{1 - x^{k + r - 1}}\prod_{\nu = 1}^{r - 1}\frac{1 - x^{k + \nu}}{1 - x^{\nu}}$$
$$= x^{k + r - 1}\left(1 + O\left(x^{k} + x^{k + r - 1}\right)\right)\prod_{\nu = 1}^{r - 1}\frac{1 - x^{k + \nu}}{1 - x^{\nu}}$$
$$ = x^{k + r - 1}\left(1 + O\left(x^{k} + x^{k + r - 1}\right) \right) G_{1}(\theta)$$
where $G_{1}(\theta)$ is closely related to the $G(\theta)$ in the generating function for $A(n,j,r)$ given in the Background Section. 

Let us now return to the standard definition of the $O$-notation.

We may apply Eq.~\ref{II}
 to the generating function $C(x)$. The product
$\prod_{\nu = 1}^{r - 1}(1 - x^{k + \nu})$ does not significantly affect the asymptotic behaviour of the derivatives of $\ln(G_{1}(\theta))$ for the same reasons that $\prod_{\nu = 1}^{r}(1 - x^{k + \nu})$ does not
significantly affect 
the derivatives of $\ln G(\theta)$. As before the asymptotic behaviour of $C(n,k,r)$ is determined by the values of its generating function at $x = e^{-\alpha}$ and its derivatives there. We find that (since $e^{-\alpha(k + r - 1)} = (\pi/\sqrt{6n})^{2}e^{-x_{1} - y_{1}}(1 + O(n^{-1/2})))$ for $k = (\sqrt{6n}/\pi)\ln(\sqrt{6n}/\pi) + y_{1}\pi^{-1}\sqrt{6n}$ and $r = (\sqrt{6n}/\pi)\ln(\sqrt{6n}/\pi) + x_{1}\pi^{-1}\sqrt{6n}$ where $|x_{1}|,|y_{1}| \le 5^{-1}\ln n$
\begin{equation}\label{VIII}
C(n,k,r) = \left(\frac{\pi}{\sqrt{6n}}\right)^{2}P(n)e^{-x_{1} - e^{-x_{1}} - y_{1} - e^{-y_
{1}}}\left(1 + O\left(n^{-1/10}\ln n\right)\right).
\end{equation}
If one of  $x_{1},y_{1}$ is $\le -5^{-1}\ln n$ then
$$C(n,k,r) = O\left(e^{-n^{1/5}}P(n)n^{-1}\right).$$
$\;$ We also have if $k \le (\sqrt{6n}/\pi)\ln(\sqrt{6n}/\pi) + y_{1}\pi^{-1}\sqrt{6n}$ and $r = \lfloor \pi^{-1}\sqrt{6n} \ln (\sqrt{6n}/\pi) + x_{1}\pi^{-1}\sqrt{6n} \rfloor$ then 
$$B(n,k,r) = \frac{\pi}{\sqrt{6n}}P(n)e^{-x_{1} - e^{-x_{1}}}e^{-e^{-y_{1}}}\left(1 + O\left(n^{-1/10}\ln n\right)\right).$$
There is a similar estimate for $B(n,k,r)$ if the restrictions on $k$ and $r$ are interchanged obtained by interchanging $x_{1}$ and $y_{1}$ in this last formula.

\section{Successive Ranks}
 The rank of a partition, first defined by Dyson\cite{Dyson}, equals the largest part minus the number of parts. The $k$-th rank, introduced by Atkin\cite{Atkin} and denoted here by $r_{k}$, equals the $k$-th largest part, $d_{k}$, minus the number of parts greater than or equal to $k$ (the $k$-th largest part in the conjugate partition), $s_{k}$. Results concerning the ranks of partitions are helpful when estimating the number of partitions of an even integer $n$ whose parts form the degree sequence of a simple(no loops or multiple edges) unlabelled graph, such partitions are are called graphical \cite{Erdos}. The degree of a vertex is the number of edges incident with it. The earliest
 result in graph theory, due to Euler, is that the sum of the vertex degrees equals twice the number of edges. Nash-Williams \cite{Nash}(unpublished) showed that if $K$ is the size(the largest $k$ such that $d_{k} \ge k$) of the Durfee square of the partition then a partition is graphical if and only if $r_{1} + r_{2} + \cdots + r_{k} \le -k$ for $k \le K$. (Canfield and Savage \cite{Canfield} show that almost all partitions of $n$ have a Durfee square of size $\sim(\sqrt{6}\ln 2/\pi)n^{1/2}$, we only consider $k$ that are $\le n^{1/10}/\ln n$ so our restriction on $k$ is irrelevant)
 See Sierksma- Hoogevan \cite{Sier} for a proof. As far as we know graphical partitions have not been enumerated. The number of partitions with $r_{1} \le -1$ has been investigated by Erd\H{o}s-Richmond \cite{Erdos}. We repeat for completeness the routine part of their analysis(in \cite{Erdos} some times the scaling $x_{1}\frac{\sqrt{6n}}{\pi}$ is used as here and sometimes the scaling 
$x_{1}\sqrt{n}$ is used as in \cite{Lehner}). In Theorem~\ref{two}
 we saw that the distributions of the largest part and the number of parts are identical independent extreme value distributions. From Feller\cite{Feller} Section 7 Eq.(5.3) the distribution of the rank is $R(t)$ where, using the substitution $v = e^{-x}$ followed by $u = v(1 + e^{-t})$,
\begin{equation}\label{IX}
R(t) = \int_{-\infty}^{\infty}e^{-e^{-t - x}}e^{-x - e^{-x}}dx = \int_{0}^{\infty}e^{-v(1 + e^{-t})}dv
\end{equation}
$$=\frac{1}{1 + e^{-t}}\int_{0}^{\infty}e^{-u}du
 = \frac{1}{1 - e^{-t}}$$
the (0,1) logistic distribution with density function $1/((1 + e^{-t})(1 + e^{t}))$. The variance of the (0,1) logistic distribution is wellknown to be $\pi^{2}/3$ since its mean is 0 and from Eq. 3 of 3.527 in 
Gradshteyn-Ryzhik \cite{Grad}
$$E(t^{2}) = \int_{-\infty}^{\infty}\frac{t^{2}}{1 + e^{t} + e^{-t} + 1}dt = 4\int_{0}^{\infty}\frac{u^{2}}{\cosh^{2}u}du = \frac{\pi^{2}}{3}.$$
This implies that the number of partitions of $n$ with rank $\le t\sqrt{6n}/\pi$ is asymptotic to $P(n)/(1 - e^{-t})$ and the number of partitions of $n$ with rank $= \lfloor t\sqrt{6n}/\pi \rfloor$ is asymptotic to 
$P(n)\frac{\pi}{\sqrt{6n}}/((1 + e^{-t})(1 + e^{t}))$ as shown in \cite{Erdos}. We should mention that the analysis in this paper has been improved by Rousseau and Ali \cite{Rousseau}.

To study the successive ranks we follow \cite{Erdos}. Let $d_{1} = \frac{\sqrt{6n}}{\pi}\ln\left(\frac{\sqrt{6n}}{\pi}\right) + x_{1}\frac{\sqrt{6n}}{\pi}$ denote the largest part of a partition $\pi$ of $n$ and let $s_{1} = \frac{\sqrt{6n}}{\pi}\ln\left(\frac{\sqrt{6n}}{\pi}\right) + y_{1}\frac{\sqrt{6n}}{\pi}$ denote the number of parts of $\pi$. We may construct a partition of $n$ with given $d_{1}$ and $s_{1}$ from a partition of $n - d_{1} - s_{1}$ iff the largest summand $d_{2}$ of the partition of $n - s_{1} - d_{1}$ satisfies $d_{2} \le d_{1}$ and similarly the number of parts $s_{2}$ of the partition of $n - d_{1} - s_{1}$ satisfies $s_{2} \le s_{1}$. Eq.~(\ref{IX}) estimates the number of partitions that have $d_{1}$ and $s_{1}$ as given in Eq.~(\ref{IX}). Now if $d_{1} + s_{1} = o(n)$ then
$$(n - d_{1} - s_{1})^{1/2} = n^{1/2}\left(1 - \frac{d_{1} + s_{1}}{n}\right)^{1/2} = n^{1/2} -\frac{d_{1} + s_{1}}{2 n^{1/2}} + O\left(\frac{(d_{1} + s_{1})^{2}}{n^{3/2}}\right)$$
$$= n^{1/2} - \frac{\sqrt{6}}{\pi}\ln\left(\frac{\sqrt{6n}}{\pi}\right) - \frac{\sqrt{3/2}}{\pi}(x_{1} + y_{1}) + O\left(n^{-3/2}(d_{1} + s_{1})^{2}\right).$$
Thus
$$\pi\sqrt{2(n - d_{1} - s_{1})/3} = \pi\sqrt{2n/3} - 2\ln\left(\pi^{-1}\sqrt{6n}\right) - x_{1} - y_{1} + O\left(n^{-3/2}(d_{1} + s_{1})^{2}\right)$$
and
$$e^{\pi\sqrt{6(n - d_{1} - s_{1})/3}} = \frac{\pi^{2}}{6n}e^{\pi\sqrt{2n/3}}e^{-x_{1}}e^{-y_{1}}\left(1 + O\left(\frac{(d_{1} + s_{1})^{2}}{n^{3/2}}\right)\right).$$
Also
$$n - d_{1} - s_{1} = n\left(1 + O\left(\frac{d_{1} + s_{1}}{n}\right)\right)$$
so if $|x_{1}|, |y_{1}| \le 5^{-1}\ln n$ we have
\begin{equation}\label{X}
P(n - d_{1} - s_{1}) = P(n)\frac{\pi^{2}}{6n}e^{-x_{1} - y_{1}}\left(1 + O\left(\frac{\ln^{2}n}{n^{1/2}}\right)\right).
\end{equation}
Consider the number of partitions of $n - d_{1} - s_{1}$ with $d_{2} \le \frac{\sqrt{6(n - d_{1} -s_{1})}}{\pi}\ln\left(\frac{\sqrt{6(n - d_{1} - s_{1})}}{\pi}\right) + x_{2}\frac{\sqrt{6(n - d_{1} - s_{1})}}{\pi}$ and
 $s_{2} \le \frac{\sqrt{6(n - d_{1} - d_{1})}}{\pi}\ln\left(\frac{\sqrt{6(n - d_{
1} - s_{1})}}{\pi}\right) + y_{2}\frac{\sqrt{6(n - d_{1} - s_{1})}}{\pi}.$
Let us examine these bounds to see how much difference the $d_{1} + s_{1}$ term makes. We will show that to the accuracy of our $O$-terms it makes no difference, can be replaced by 0 that is. If
$$f(x) = \frac{(6(n - x))^{1/2}}{\pi}\ln(\pi^{-1}(6(n - x))^{1/2}) + (x_{2} + y_{2})\pi^{-1}\sqrt{6(n - x)}$$
then for $|x_{2}|, |y_{2}| \le 5^{-1}\ln n$
$$f^{'}(0) = -\frac{\sqrt{\frac{3}{2n}}}{\pi}
\left(\ln\left(\frac{\sqrt{6n}}{\pi}\right) + 1 + x_{2} + y_{2}\right) = O(n^{-1/2}\ln n)$$
so
$$f(d_{1} + s_{1}) = \frac{\sqrt{6n}}{\pi}\ln\left(\frac{\sqrt{6n}}{\pi}\right) + \left(x_{2} + y_{2} + O(n^{-1}\ln n)\right)\frac{\sqrt{6n}}{\pi}.$$
We can view the $O(n^{-1}\ln n)$ term as a change in $x_{2}$ and or $y_{2}$. It produces a $O(n^{-1}\ln n)$ change in the functions
$e^{-e^{-x_{2}}}$ and $e^{-e^{-y_{2}}}$. Thus the number of partitions of $n - d_{1} - s_{1}$ with $d_{2} \le \frac{\sqrt{6n}}{\pi}\ln\left(\frac{\sqrt{6n}}{\pi}\right)
+ x_{2}\frac{\sqrt{6n}}{\pi}$ and $s_{2} \le \frac{\sqrt{6n}}{\pi}\ln(\frac{\sqrt{6n}}{\pi}) +
 y_{2}\frac{\sqrt{6n}}{\pi}$ for $x_{1} \ge x_{2}, y_{1} \ge y_{2}$ from Eq.\ref{X}
 and Theorem~\ref{one} for $|x_{1}|, |y_{1}| \le 5^{-1}\ln n$ and $x_{1} \ge x_{2}, y_{1} \ge y_{2}$ is asymptotic to
$$P(n - d_{1} - s_{1})e^{-e^{-x_{2}}}e^{-e^{-y_{2}}}\left(1 + O\left(\frac{\ln^{2}n}{n^{1/10}}\right)\right)$$
$$ = P(n)e^{-x_{1} -y_{1} - e^{-x_{2}} - e^{-y_{2}}}
\left(1 + O\left(\frac{\ln^{2}n}{n^{1/10}}\right)\right).$$
Of course as we saw in the proof of Theorem~\ref{two} if one of $x_{1},x_{2}, y_{1},y_{2}$ is $\le -5^{-1}\ln n$ then the righthand side of this last equation is exponentially small compared to $P(n)$. This equation agrees with the result of Fristedt \cite{Fristedt} although Fristed's result is much more general however Fristedt does not consider the independence of the distributions defined by the the $d$-th largest part in a partition and its conjugate.
In the proof of his Theorem 2.5 Fristedt integrates, so we do not repeat the integration here however we may do the same integrations for the independent distributions defined by the $k$-th largest part in a partition and its conjugate as Fristedt does, 
$$exp\left(- e^{-v_{t}} - \sum_{s = 1}^{t - 1}v_{s}\right) \;\; \mbox{for}\;\; v_{1}\ge v_{2} \ge \cdots \ge v_{t}$$
over all possible values and obtains his Theorem 2.3 where $Y_{d}$ denotes the distribution defined by the $d$-th largest part, which states
\begin{equation}\label{XI}
Pr\left(\frac{\pi}{\sqrt{6n}}Y_{d} - \ln(\pi^{-1} \sqrt{6n}) \le y \right) \rightarrow \int_{-\infty}^{y}\frac{e^{-e^{-v} -dv}}{(d - 1)!}dv.
\end{equation}
 Fristedt shows that with high probability, ie. tending to 1 with $n$, that the largest parts are distinct. 
Since if $d = o(n^{1/10}\ln^{-2}n)$ the distributions of the $d$-th largest part and the number of parts $\ge d$ are equal and independent this will also be true for the largest parts in the conjugate partitions for this range of $d$.

Note Eq. ~(\ref{X}) can be easily generalized to the partitions of $n - d_{1} - d_{2} - \cdots -d_{t} - s_{1} - s_{2} - \cdots -s_{t}$ so we find from Fristedt's results and our methods that for $d = o(n^{1/10}\ln^{-2}n)$ the number of partitions with $d_{i} = \frac{\sqrt{6n}}{\pi} + x_{i}\frac{\sqrt{6n}}{\pi}, s_{i} = \frac{\sqrt{6n}}{\pi} + y_{i}\frac{\sqrt{6n}}{\pi}, 1 \le i \le t \le d,$ is asymptotic to 
$$P(n)e^{\sum_{i = 1}^{t-1} - x_{i} - y_{i} - e^{-x_{t}} - e^{-y_{t}}}.$$
Note if $|f_{i}| \le n^{-1/10}\ln^{-2}n$ then $\prod_{i = 1}^{K}(1 + f_{i}) = \sum_{i = 1}^{K}f_{i} + O\left(\sum_{i = 1}^{K}f_{i}^{2}\right).$
Furthermore we have that with $Y_{d}(x)$ denoting the distribution in Eq.~(\ref{XI})(Eq.~(\ref{XI}) shows that the $Y_{d}$ are not identical)
\begin{thm}\label{three}
Let the $i$-th largest part of a partition of $n$ be bounded by $d_{i} \le \frac{\sqrt{6n}}{\pi}\ln\left(\frac{\sqrt{6n}}{\pi}\right) + x_{i}\frac{\sqrt{6n}}{\pi}$ and if the number of parts $\ge i$

is bounded by $s_{i} \le \frac{\sqrt{6n}}{\pi}\ln\left(\frac{\sqrt{6n}}{\pi}\right) + y_{i}\frac{\sqrt{6n}}{\pi}$ then the number of partitions satisfying these bounds for $i \le K$, say $P(n, \le x_{i}, \le y_{i})$ satisfies
$$P(n, \le x_{i}, \le y_{i}) = P(n)\prod_{i = 1}^{K}Y_{i}(x_{i})Y_{i}(y_{i})\left(1 + O\left(\frac{K\ln^{2}n}{n^{1/10}}\right)\right).$$
Again the distributions defined by $Y_{i}(x_{i})$ and $Y_{i}(y_{i})$ are independent and equal for each $i$. If $i \ne j$ then $Y_{i}(x_{i}), Y_{j}(y_{j})$ are not identical but are independent.
\end{thm}

The probability that one partition of $n$ dominates another is with $\lambda_{i}$ denoting $d_{i}$ in one partition $\pi$ and $\lambda_{i}^{'}$ denoting the $i$-th part in the conjugate partition of $\pi$ is the probability that $\sum_{i =1}^{k}\lambda_{i} - \lambda_{i}^{'} \le 0$(Nash-Williams' condition is that this sum is $\le -k$ for each $k \le d(\lambda)$, the size of the Durfee square of $\pi$). We cannot say anything new about Macdonald's conjecture that the probability that one partition of $n$ dominates another tends to 0 however Esseen's result in probability theory allows us to say more about Nash-William's 
 criterion. H. Wilf conjectured that the fraction of the partitions of $n$ that are graphical tends to zero with $n$ and Pittel \cite{Pittel} verifies this conjecture and Macdonald's
but without a rate of convergence estimate. We have seen that the distribution defined by the first rank is a logistic (0,1) distribution. We now state Esseen's result in Feller \cite{Feller} Theorem 2, p. 544.

Let the $X_{k}$ be independent variables such that
$$E(X_{k}) = 0, \;\; E(X_{k}^{2}) = \sigma_{k}^{2} > 0,\;\; E(|X_{k}|^{3}) = \rho_{k} < \infty.$$
Put $s_{n}^{2} = \sigma_{1}^{2} + \sigma_{2}^{2} + \cdots + \sigma_{n}^{2}, \;\; r_{n} = \rho_{1} + \rho_{2} + \cdots + \rho_{n}$
and let $F_{n}$ stand for the distribution of the normalized sum
$$(X_{1} + \cdots + X_{n})/s_{n}.$$
Then for all $x$ and $n$
$$|F_{n}(x) - \mathbf{N}(x)| \le 6\frac{r_{n}}{s_{n}^{3}}.$$

Now from Ch. 7 of Abromowitz and Stegun \cite{Abrom}, Eq. 7.1.23 and Eq. 7.1.2
$$\mathbf{N}(x) = \frac{1}{\sqrt{2\pi}}\int_{-\infty}^{x}e^{-t^{2}/2}dt = \frac{1}{\sqrt{\pi}}\int_{x/\sqrt{2}}^{\infty}e^{-u^{2}}du = \frac{1}{2}erfc(x/\sqrt{2})$$ 
$$= \frac{e^{-x^{2}/2}}{x\sqrt{2\pi}}\left(1 + O\left(x^{-2}\right)\right).$$
Now from Abromovitz and Stegun, \cite{Abrom} 
$$\int_{-\infty}^{-k}e^{-t^{2}/2}dt = \int_{k}^{\infty}e^{-t^{2}/2}dt = \sqrt{2}\int_{k/\sqrt{2}}^{\infty}e^{-u^{2}}du $$
$$=\frac{\sqrt{\pi}}{\sqrt{2}}erfc(k/\sqrt{2}).$$
From Eq. 7.1.23 of \cite{Abrom}
$$erfc(z) = \frac{1}{z\sqrt{\pi}}e^{-z^{2}}\left(1 + O\left(z^{-2}\right)\right).$$ 

 We shall find it convenient to have a finite explicit formula for the $Y_{d}$ of Eq.~(\ref{XI}). Such a formula follows by integrating by parts. If we again let $u = e^{-v}$ in Eq.~(\ref{XI}) we find that 
$$Y_{k}(x) = \frac{1}{\Gamma(k)}\int_{-\infty}^{x}e^{-kv - e^{-v}}dv = \frac{1}{\Gamma(k)}\int_{e^{-x}}^{\infty}u^{k - 1}e^{-u}du$$
$$\left. = \frac{1}{\Gamma(k)}u^{k - 1}\frac{d-e^{-u}}{du}\right|_{e^{-x}}^{\infty}
+\frac{1}{\Gamma(k - 1)} \int_{e^{-x}}^{\infty}u^{k - 2}e^{-u}du.$$
If we continue integrating by parts we find that
\begin{equation}\label{XII}
Y_{k}(x) = e^{-e^{-x}}\sum_{i = 1}^{k}\frac{e^{-(k - i)x}}{\Gamma(k - i +1)}
\end{equation}
Since the distributions defined by $d_{k}$ and $s_{k}$, let us say $d_{k}(x)$ and $s_{k}(x)$ are defined by Eq.~(\ref{XI}) we have from Eq. (5.3) of Section V.5 of Feller \cite{Feller} that the distribution of
$d_{k}(x) - s_{k}(x) = r_{k}(x)$ is defined by, recall that $Y_{k}^{'}(x) = e^{-kx - e^{-x}}/\Gamma(k)$,
$$r_{k}(t) = \int_{-\infty}^{\infty}e^{-e^{-x -t }}\sum_{i = 1}^{k}\frac{e^{-(k - i)x - (k - i)t}}{\Gamma(k - i + 1)} 
\frac{e^{-kx - e^{-x}}}{\Gamma(k)}dx$$
$$ = \frac{1}{\Gamma(k)}\int_{-\infty}^{\infty}e^{e^{-x}(1 + e^{-t})}\sum_{i = 1}^{k}\frac{e^{-(2k - i)x}}{\Gamma(k - i + 1)}dx.$$
If we again let $u = e^{-x}, dx = -du/u$ we find
$$r_{k}(t) = \frac{1}{\Gamma(k)}\int_{0}^{\infty}e^{-u(1 + e^{-t})}\sum_{i = 1}^{k}\frac{u^{2k - i - 1}}{\Gamma(k - i + 1)}e^{-(k - i)t}du.$$
Again let $v = (1 + e^{-t})u, du = (1 + e^{-t})^{-1}dv$ and we find that
$$r_{k}(t) = \frac{1}{\Gamma(k)\left(1 + e^{-t}\right)^{2k}} 
\int_{0}^{\infty}e^{-v}\sum_{i = 1}^{k}\frac{v^{2k - i - 1}}{\Gamma(k - i + 1)}e^{-(k - i)t}(1 + e^{-t})^{i}dv$$
$$ = \frac{1}{\Gamma(k)\left(1 + e^{-t}\right)^{2k}}\sum_{i = 1}^{k}\frac{\Gamma(2k - i)}
{\Gamma(k - i + 1)}\left(1 + e^{-t}\right)^{i}e^{-(k - i)t}.$$
To find a nice formula for the derivative of $r_{k}(t)$ we write this as 
$$ r_{k}(t) = \frac{1}{\Gamma(k)} \sum_{i = 1}^{k}\frac{\Gamma(k + i -1)e^{-(i - 1)t}}{\Gamma(i)\left(1 + e^{-t}\right)^{k + i - 1}}.$$
Now
$$r_{k}^{'}(t) = \frac{ke^{-t}}{\left(1 + e^{-t}\right)^{k + 1}} - \frac{ke^{-t}}{\left(1 + e^{-t}\right)^{k + 1}} + \frac{k(k + 1)e^{-2t}}{\Gamma(2)\left(1 + e^{-t}\right)^{k + 2}}$$
$$ - \frac{k(k + 1)e^{-2t}}{\Gamma(2)\left(1 + e^{-t}\right)^{k + 2}} + \cdots + \frac{\Gamma(2k -1)(2k - 1)e^{-kt}}{\Gamma^{2}(k)\left(1 + e^{-t}\right)^{2k}}$$
The terms telescope and we conclude that
$$r_{k}^{'}(t) = \frac{\Gamma(2k)e^{-kt}}{\Gamma^{2}(k)\left(1 + e^{-t}\right)^{2k}}.$$
This shows that the number of partitions of $n$ with $r_{k}(t) = \lfloor \frac{\sqrt{6n}}{\pi}t_{k}\rfloor$ is
$$\sim P(n)\frac{\pi}{\sqrt{6n}}\frac{\Gamma(2k)}{\Gamma^{2}(k)}\frac{1}{\left(1 + e^{-t_{k}}\right)^{k}\left(1 + e^{t_{k}}\right)^{k}}.$$
This generalizes the results for $r_{1}$ in \cite{Erdos} and moreover since each pair $d_{k}(x)$ and $s_{k}(x)$ is
independent of all other pairs $(d_{j}(x),s_{j}(x)), j \ne k$ we see that the distributions $r_{k}(x)$ are 
independent so we may generalize our asymptotic formula for the number of partitions of $n$ with $r_{k}(t) = \lfloor \frac{\sqrt{6n}}{\pi} t_{k}\rfloor$ to an asymptotic formula for the number of partitions of $n$ with any specified $r_{k}(t)$ for any set of $k$ that satisfy that $k = o(n^{1/10}\ln^{-2}n)$. Perhaps the result is too cumbersome to state however.

 To use Esseen's result we calculate the variance and the expectation of $|t|^{3}$.
Since from the asymptotic formula for $\Gamma(k)$ we have
$$r_{k}^{'}(t) = 4^{k}k^{1/2}\sqrt{\pi}\left(\frac{e^{-t}}{\left(1 + e^{-t}\right)^{2}}\right)^{k}\left(1 + O\left(\frac{1}{k}\right)\right)$$
we consider
$$\frac{e^{-t}}{\left(1 + e^{-t}\right)^{2}} = exp(-t -2\ln\left(1 + e^{-t}\right)) = exp\left(-t -2\ln\left(2 - t + t^{2}/2 + \cdots \right)\right)$$
$$= exp\left(-2\ln 2 -t^{2} + O(t^{3}\right) = \frac{e^{-t^{2}}}{4}\left(1 + O\left(t^{3}\right)\right).$$

Thus
\begin{equation}\label{XIII}
\left(\frac{e^{-t}}{\left(1 + e^{-t}\right)^{2}}\right)^{k} = \frac{e^{-k t^{2}}}{4^{k}}\left(1 + O\left(kt^{3}\right)\right).
\end{equation}
We have
$$\sigma_{k}^{2} = \frac{\Gamma(2k)}{\Gamma^{2}(k)}2\int_{0}^{\infty}\frac{t^{2}}{\left(1 + e^{-t}\right)^{k}\left(1 + e^{t}\right)^{k}}dt.$$
Let us consider the integrals over $[0, k^{-2/5}]$ and $[k^{-2/5}, \infty]$ separately.
 $$\int_{k^{-2/5}}^{\infty}t^{2}\left(\frac{e^{-t}}{\left(1 + e^{-t}
\right)^{2}}\right)^{k}dt < \int_{k^{-2/5}}^{\infty}t^{2}e^{-kt}dt$$
$$ =  \frac{1}{k^{3}}\int_{k^{3/5}}^{\infty}u^{2}e^{-u}du = O\left(k^{-9/5}e^{-k^{3/5}}\right).$$
Also from Eq.~(\ref{XIII})
$$2\int_{0}^{k^{-2/5}}t^{2}\left(1 + e^{-t}\right)^{-k}\left(1 + e^{t}\right)^{-k}dt = 2\times 4^{-k}\int_{0}^{k^{-2/5}}t^{2}e^{-kt^{2}}\left(1 + O\left(kt^{3}\right)\right)$$
$$= k^{-3/2}4^{-k}2\int_{0}^{k^{1/10}}u^{2}e^{-u^{2}}du = 2\times 4^{-k}k^{-3/2}\left(\frac{\sqrt{\pi}}{4}
- \int_{k^{1/10}}^{\infty}u^{2}e^{-u}du\right)$$
$$ = 4^{-k}k^{-3/2}\sqrt{\pi}/2 + O\left(k^{-13/10}e^{-k^{1/5}}\right).$$
Thus $\sigma^{2}_{k} \sim \pi k^{-1}/2$ and the $s_{n}^{2}$ in Esseen's result is $\sim \pi 2^{-1} \ln n$.
The asymptotic behaviour of the $\rho_{k}$ and $r_{n}$ in Esseen's result may be determined in much the same way, the $t^{2}$ in our integrals for $\sigma^{2}_{n}$ are replaced by $t^{3}$. The result is that $\rho_{k}$ is asymptotic to $12 k^{-3/2}\sqrt{\pi}$. Thus $r_{n} \sim 12\sqrt{\pi}\zeta(3/2).$ It follows that $r_{n}s_{n}^{-3/2} \sim 24\sqrt{2}\pi^{-1}\zeta(3/2)\ln^{-3/2}n$.
Since the integral from $-\infty$ to $-k$ of the density of the normal distribution is exponentially small we deduce from Esseen's result that the probability that a partition is graphical is $O(\ln^{-1/2}n)$ since $Pr(\sum_{l = 1}^{k}r_{l} \le -k)s_{n}^{-1} = O(s_{n}^{-3/2})$ if and only if $Pr(\sum_{l = 1}^{k}r_{l} \le -k) = O(s_{n}^{-1/2}) = O(\ln^{-1/2} n)$ and we have 
\begin{thm}
The probability that the parts of a partition of an even integer $n$ are the degree sequence of a simple graph is
$$O\left(\frac{1}{\ln^{1/2}n}\right).$$
\end{thm}

\section{Acknowledgements and Comments}
We are greatly indebted to George Andrews for the Tak\'{a}cs reference and helpful comments on presentation. Of course any errors and faults in presentation are solely the author's responsibility. It would be interesting to have a generating function for graphical partitions. Perhaps this would allow the determination of an asymptotic formula for the number of graphical partitions of $n$. Probably our upper bound is too big. The best known lower bound is due to Rousseau and Ali \cite{Rousseau} and is of the form $cP(n)n^{-1/2}$. Since we end up estimating the probability that a sum of distributions is in the tail of the sum distribution it is likely our upper bound can be improved by someone knowing probability theory?


\begin{thebibliography}{99}
\bibitem{Abrom}

M. Abromowitz and I. A. Stegun,
\emph{Handbook of Mathematical Functions}, Dover, 1970.

\bibitem{AA}
Gert Almkvist and George E. Andrews,
\emph{A Hardy-Ramanujan Formula for Restricted Partitions}, J. Number Theory, 38, 135-144 (1991).

\bibitem{Andrews}
G. E. Andrews,
\emph{The Theory of Partitions}, vol. 2 of \emph{Encyclopedia of Mathematics and its Applications}, Addison-Wesley, 1976.

\bibitem{Atkin}
A. O. L. Atkin,
\emph{A note on ranks and conjugacy of partitions}, Quart. J. Math. Oxford Ser(2), 17 (1966), 335-338.

\bibitem{Canfield}
E. R. Canfield and C. Savage,
\emph{Durfee Polynomials}, The Electronic J. of Combinatorics, 5 (1998).

\bibitem{Comtet}
Louis Comtet,
\emph{Advanced Combinatotics}, translated fom the French by J. W. Neinhuys, D. Reidel, 1974.

\bibitem{Dyson}
F. J. Dyson, \emph{Some guesses in the theory of partitions}, Eurika (Cambridge), 8 (1944), 10-15.

\bibitem{Lehner}
P. Erd\H{o}s and J. Lehner,
\emph{The distribution of the number of summands in the partitions of a positive integer}, Duke Math. J. 8(1941), 335-45.

\bibitem{Erdos}
P. Erd\H{o}s and L. B. Richmond,
\emph{On Graphical Partitions}, Combinatorica, 13(1), 57-63, (1993).

\bibitem{Feller}
William Feller,
\emph{An Introduction to Probability Theory and Its Applications}, 2nd ed., J. Wiley and Sons, 1971

\bibitem{Flajolet}
Philippe Flajolet and Robert Sedgewick,
\emph{Analytic Combinatorics}, Cambridge University Press, 2009.

\bibitem{Fristedt}
B. Fristedt, \emph{The structure of random partitions of large integers}, Trans. Amer. Math. Soc. 337 (1993), 703-735.

\bibitem{Gauss}
Gauss, C. F. , \emph{Werke}, vol. 2, K\"{o}nigliche Gesellschaft der Wissenschaften, Gottingen, 1863.

\bibitem{Gould}
I. P. Goulden and D. M. Jackson, \emph{Combinatorial Enumeration}, J. Wiley and Sons, 1983.

\bibitem{Grad}
I. S. Gradshteyn and I. M. Ryzhik,
\emph{Table of Integrals, Series and Products}, 4th ed., Academic Press, 1965.

\bibitem{Macdonald}
I. G. Macdonald,
\emph{Symmetric Functions and Hall Polynomials}, 2nd edition, Oxford Sci. Publ. Oxford Univ. Press, London OX2 6DP, 1995.

\bibitem{Nash}
C. St. J. A. Nash-Williams: Valency sequences which force graphs to have Hamiltonian circiuts; Interim Report C. and O. Research Report, Fac. of Math., U. of Waterloo.

\bibitem{Pittel}
B. Pittel, \emph{Confirming Two Conjectures about the Integer Partitions}, JCT(A) 88, 123-135 (1999).

\bibitem{Rousseau}
C. Rousseau and F. Ali, \emph{On a conjecture concerning graphical partitions}, Congr. Num. 104(1994), 150-160.

\bibitem{Stanley}
R. Stanley, \emph{Enumerative Combinatorics, Vol. 1}, Cambridge University Press, 1997.

\bibitem{Sier}
G. Sierksma and H. Hoogeveen, \emph{Seven criteria for integer sequences being graphic}, J. Graph Theory, 15(1991), p. 96-111.

\bibitem{Sz1}
George Szekeres,
\emph{An asymptotic formula in the theory of partitions, II}, Quart. J. Math. Oxford Se. (2), 4.(1951), p. 96-111.

\bibitem{Sz2}
George Szekeres,
\emph{Asymptotic distribution of the number and size of parts in unequal partitions}, Bull. Austral. Math. Soc., vol. 36 (1987) 89-97.

\bibitem{Takacs}
Lajos Tak\'{a}cs, \emph{Some Asymptotic Formulas for Lattice Paths}, J. Stat. Plann. and Inference 14 (1986), 123-142.

\end{thebibliography}
\end{document}